\documentclass{revtex4}
\usepackage{graphicx}

\usepackage{amsfonts}
\usepackage{amsmath}

\begin{document}


\title{Noise synchronisation and stochastic bifurcations in lasers}

\author{Sebastian M. Wieczorek}

\affiliation{Mathematics Research Institute,
   University of Exeter, Exeter EX4 4QF, UK, s.m.wieczorek@ex.ac.uk}

\maketitle

\section{Introduction}
\label{sec:intro}

Synchronisation of nonlinear oscillators to  irregular external signals is an interesting problem
of importance in physics, biology, applied science, and 
engineering~\cite{PIK01,PIK84,MAI95,JEN98,ZHO03,JUN04,UCH04,KOS04,GOL05,GOL06}.  
The key difference to synchronisation to a periodic external signal 
is the lack of a simple functional relationship between the  input signal and the synchronised 
output signal, making the phenomenon much less evident~\cite{PIK01,RUL95}. Rather, 
synchronisation is detected when  two or more identical uncoupled oscillators 
 driven by the same external signal but starting at 
different initial states  have identical long-term responses. This is equivalent to obtaining 
reproducible long-term response from a single oscillator driven repeatedly by the same 
external signal, each time starting at a different initial state. Hence, synchronisation to
irregular external signals is also known as reliability~\cite{MAI95} or 
consistency~\cite{UCH04}, and represents the ability to encode irregular signals in a 
reproducible manner.

Recent studies have shown that nonlinear oscillators can exhibit interesting  responses 
to stochastic external signals.  Typically, a small amount of external noise causes 
synchronisation~\cite[Ch.7]{PIK01},~\cite{PIK84,JEN98}. However, as the strength of external noise 
increases, there can be a loss of synchrony in oscillators with amplitude-phase coupling 
(also known as  shear, nonisochronicity or amplitude-dependent frequency)~\cite[Ch.7]{PIK01},
~\cite{KOS04,GOL05,GOL06,WIE09,LIN07,LIN08}. 
Mathematically, loss of noise synchrony, consistency or reliability is a manifestation of
a {\it stochastic bifurcation of a random attractor}.

This chapter gives a definition of noise synchronisation in terms of  random 
pullback attractors and studies synchronisation-desynchronisation transitions as purely 
noise-induced stochastic 
bifurcations. This is in contrast to the effects described in Chapter 2, where noise is used to 
control or regulate the dynamics that is already present in the noise-free system. 
We focus on a single-mode class-B laser model  and the Landau-Stuart model 
(Hopf normal form with shear~\cite{PIK01}). 
In section~\ref{sec:bifdiag}, numerical analysis of  the locus 
of the stochastic bifurcation in a three-dimensional parameter space of  the `distance' 
from Hopf  bifurcation,  amount of amplitude-phase coupling, and  external signal strength reveals a 
simple power law for the Landau-Stuart model but quite different behaviour for the laser model.
In section~\ref{sec:interpret}, the analysis  of the shear-induced stretch-and-fold action 
that creates horseshoes
gives an intuitive explanation for the observed loss of synchrony and for the deviation from the 
simple power law in the laser model. Experimentally, stochastic external forcing can be realised 
by optically 
injecting noisy light into a (semiconductor) laser as described in section~\ref{sec:synch}.
While bifurcations of random pullback attractors and the associated synchronisation-desynchronisation 
transitions have been studied theoretically, single-mode semiconductor lasers emerge as interesting 
candidates for experimental testing of these phenomena.

\section{Class-B laser model and Landau-Stuart model}
\label{sec:model}

A class-B single-mode laser~\cite{WIE05} without noise can be modelled by the rate equations~\cite{WIE08}:
\begin{eqnarray}
\label{eq:E}
\frac{dE}{dt}&=& i\Delta\,E+ g\gamma(1-i\alpha)NE,\\
\label{eq:N}
\frac{dN}{dt}&=& J - N - (1+g N)|E|^2,
\end{eqnarray}
which define a three-dimensional dynamical system with a normalised  electric field amplitude, 
$E\in{\mathbb C}$, and normalised deviation from the threshold population inversion, $N\in{\mathbb R}$, 
such that $N=-1$ corresponds to zero population inversion. 
Parameter $J$ is the normalised deviation from the threshold pump rate such that $J=-1$ 
corresponds to zero pump rate. The linewidth enhancement 
factor, $\alpha$, quantifies the amount of amplitude-phase coupling, $\Delta$ is 
the normalised detuning (difference) between  some conveniently 
chosen reference frequency and the natural laser frequency,  $\gamma=500$ is the 
normalised decay rate, and $g=2.765$ is the normalised gain coefficient~\cite{WIE08}.

System~(\ref{eq:E}--\ref{eq:N})   is $\mathbb{S}$$^1$-equivariant, meaning that it has 
rotational  symmetry
corresponding to a phase shift $E\rightarrow Ee^{i\phi}$, where $0<\phi\le 2\pi$. 
For $J\in\mathbb{R}$, there is an equilibrium at $(E,N)=(0,J)$ 
which represents the ``off'' state of the laser. This equilibrium is globally stable 
if $J<0$ and unstable if $J>0$.  At $J=0$, there is a Hopf ($\Delta\ne 0$) 
or pitchfork ($\Delta = 0$) bifurcation defining the laser threshold. 
Moreover, if $J>0$, the system has a stable group orbit in the form of
periodic orbit for $\Delta\ne 0$ or a circle of infinitely many non-hyperbolic 
(neutrally stable) equilibria for $\Delta=0$. 
In this paper, we refer to this circular attractor as the {\it limit cycle}. The limit cycle is given by 
$(|E|^2,N)=(J,0)$ and represents the ``on'' state of the laser.
Owing to the $\mathbb{S}$$^1$-symmetry,  the Floquet exponents of the limit cycle can be calculated 
analytically as eigenvalues of one of the non-hyperbolic equilibria for $\Delta=0$. 
Specifically, if 
$$
0<J<\frac{4\gamma\left(1-\sqrt{1-1/(2\gamma)}\right) -1}{g}\approx 9\times10^{-5},
$$
the overdamped limit cycle  has three real Floquet exponents
\begin{eqnarray}
\label{eq:mu1}
\mu_1=0,\;\;\;\mu_{2,3}=-a\pm b,
\end{eqnarray}
and if 
$$
\frac{4\gamma\left(1-\sqrt{1-1/(2\gamma)}\right) -1}{g}<J<\frac{4\gamma\left(1+\sqrt{1-1/(2\gamma)}\right) -1}{g}\approx1446,
$$
the underdamped limit cycle has one real and two complex-conjugate Floquet exponents
\begin{eqnarray}
\label{eq:mu2}
\mu_1=0,\;\;\; \mu_{2,3}=-a\pm ib,
\end{eqnarray}
where 
$$
a=\frac{1}{2}\left(1+g J\right)>0\;\;\;\;\;\mbox{and}\;\;\;\;\; b=\sqrt{|a^2-2g\gamma J|}>0.
$$
In the laser literature, the decaying oscillations found for pump rate in the realistic range 
$J\in(9\times 10^{-5},20 )$ are called {\it relaxation oscillations}. 
(This should not be confused with a different phenomenon of  self-sustained, slow-fast oscillations.)
Finally, even though the laser model~(\ref{eq:E}--\ref{eq:N}) is three dimensional, it cannot admit chaotic solutions
due to restrictions imposed by the rotational symmetry.

Using centre manifold theory~\cite{GUC83}, the dynamics 
of~(\ref{eq:E}--\ref{eq:N})  near the Hopf bifurcation can be approximated 
by  the two-dimensional invariant centre manifold
$$
W^c=\{(E,N)\in\mathbb{R}^3: N=J-|E|^2\},
$$
on which~(\ref{eq:E}--\ref{eq:N}) reduces to 
\begin{eqnarray}
\label{eq:null}
\frac{1}{g\gamma}\frac{dE}{dt}= 
\left[J +i\left(\frac{\Delta}{g\gamma}-\alpha(J -|E|^2) \right) \right]E - E|E|^2\nonumber.
\end{eqnarray}
After rescaling time  and detuning,
$$
\tilde{t}=t g\gamma\;\;\;\;\;\mbox{and}\;\;\;\;\; \tilde{\Delta}=\Delta/(g\gamma),
$$
we obtain the Landau-Stuart model
\begin{eqnarray}
\label{eq:reduced}
\frac{dE}{d\tilde{t}}= 
\left[J +i\left(\tilde{\Delta}-\alpha(J -|E|^2) \right) \right]E - E|E|^2,
\end{eqnarray}
that is identical to the Hopf normal form~\cite{GUC83} except for the higher order term, 
$\;i\alpha(J -|E|^2)E$, representing amplitude-phase coupling. Since this term 
does not affect stability properties of~(\ref{eq:reduced}), it does not appear in the 
Hopf normal form. However, in the presence of an external forcing, $f_{ext}(t)$,
this term has to be included because it gives rise to qualitatively different dynamics
for different values of $\alpha$. If $J>0$, the  Landau-Stuart model has a stable limit cycle with
two Floquet exponents
\begin{eqnarray}
\label{eq:fels}
\mu_1=0\;\;\;\;\mbox{and}\;\;\;\;\mu_2=-2J.
\end{eqnarray}


\begin{figure}[t]
\includegraphics[width=11.cm]{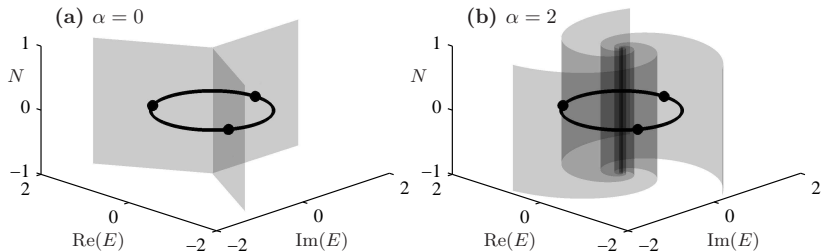}
\caption{
(Black) The limit cycle  representing the "on" state of the laser for $J=1$ and 
(gray) isochrones for three different points on the limit cycle as
defined by Eq.~(\protect\ref{eq:isochrone}) for  (a) $\alpha=0$  and  (b) $\alpha=2$.
}
\label{fig:1}
\end{figure}

\section{The linewidth enhancement factor and shear}

The linewidth enhancement factor, $\alpha$, quantifying the 
amount of amplitude-phase coupling for the complex-valued electric 
field, $E$, is absolutely crucial to our analysis.
Its physical origin is the dependence of the semiconductor refractive index, 
and hence the laser-cavity resonant frequency, on the population inversion~\cite{HEN83,WIE05}. 
A change in the electric field intensity, 
$\delta|E|^2$, induces a change,  $\delta N$, in population inversion
[Eq.~(\ref{eq:N})]. The resulting change in the refractive index  
shifts the cavity resonant frequency. The ultimate result  is a change of  $-\alpha g \gamma \delta N$ in
the {\it instantaneous frequency} of the electric field defined as $d(\arg(E))/dt$. 

Mathematically, amplitude-phase coupling 
is best illustrated by an invariant set associated with each point, $q$, on the limit 
cycle. For a point $q(0)$ on a stable limit cycle in a $n$-dimensional system, this set is defined as
\begin{eqnarray}
\label{eq:isochrone}
\{ x(0)\in\mathbb{R}^n: x(t)\rightarrow q(t)\;\;\mbox{as}\;\;t \rightarrow \infty \},
\end{eqnarray}
and is called an {\it isochrone}~\cite{GUC74}. In the laser model~(\ref{eq:E}--\ref{eq:N})
and Landau-Stuart model~(\ref{eq:reduced}), isochrones are 
logarithmic spirals that satisfy
\begin{eqnarray}
\label{eq:isochrone3}
\arg(E) +\alpha\,\ln|E|=C,\;\;\;\;\mbox{where}\;\;\;\; C\in(0,2\pi].
\end{eqnarray}
To see this, define a phase 
\begin{eqnarray}
\label{eq:isochrone2}
\Psi=\arg(E) +\alpha\,\ln|E|,
\end{eqnarray}
and check that $d\Psi/dt$ is constant and equal to $\Delta$ for~(\ref{eq:E}--\ref{eq:N}) and $\tilde{\Delta}$
for~(\ref{eq:reduced}). This means that trajectories for different initial conditions with identical initial phase, $\Psi(0)$,
will retain identical phase, $\Psi(t)$, for all time $t$. Since the limit cycle is stable, all such trajectories will converge 
to the limit cycle, where they have the same $|E(t)|$. Then,
Eq.~(\ref{eq:isochrone2}) implicates that all such trajectories have the same $\arg(E(t))$ and hence converge to just one special trajectory along the 
limit cycle as required by~(\ref{eq:isochrone}).

Isochrones of three different points on the laser limit cycle 
are shown in Fig.~\ref{fig:1}. 
Isochrone inclination to the direction normal to the limit cycle at $q(0)$  
indicates the strength of phase space stretching along the limit cycle. 
If $\alpha=0$,  trajectories with different $|E|>0$ rotate  
around the origin of the $E$-plane   with the same angular frequency giving 
no isochrone  inclination and hence no phase space stretching [Fig.~\ref{fig:1}(a)]. 
However, if $|\alpha|>0$, trajectories with larger  $|E|$ rotate
with higher angular frequency  giving rise to 
isochrone inclination and phase space stretching [Fig.~\ref{fig:1}(b)]. 
Henceforth, we refer to amplitude-phase coupling as {\it shear}.

\begin{figure}[t!]
\includegraphics[width=11cm]{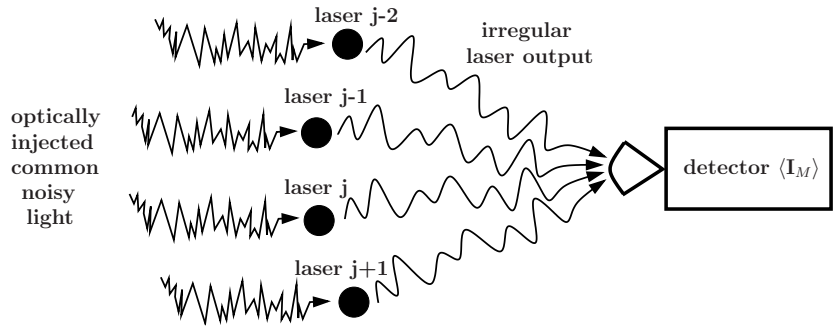}
\caption{
An experimental setup for detecting noise synchronisation  in lasers.
}
\label{fig:1b}

\end{figure}

\section{Detection of Noise Synchronisation}
\label{sec:synch}

There are at least two approaches to detecting synchronisation of a semiconductor laser to an 
irregular external signal. One approach involves a comparison of the responses of two or more 
identical and uncoupled lasers that are driven by the same external signal. The other 
approach involves a comparison of the responses of a single laser driven repeatedly by the 
same external signal~\cite{UCH04}. Here, we consider responses of $M$ uncoupled lasers 
with intrinsic spontaneous emission noise 
that are subjected to {\it common} {\it optical external forcing}, $f_{ext}(t)$,~\cite{WIE08,LAN82}:
\begin{eqnarray}
\label{eq:E2}
\frac{dE_j}{dt}&=& i\Delta\,E_j+ g\gamma(1-i\alpha)N_jE_j + f_{Ej}(t) + f_{ext}(t)\\
\label{eq:N2}
\frac{dN_j}{dt}&=& J - N_j - (1+g N_j)|E_j|^2 + f_{Nj}(t),\\
j&=&1,2,\ldots,M,\nonumber
\end{eqnarray}
The lasers are identical except for the intrinsic spontaneous emission
noise that is represented by random Gaussian processes 
$$
f_{Ej}(t)= f_{Ej}^R(t)+if_{Ej}^I(t)\;\;\;\;\mbox{ and}\;\;\;\; f_{Nj}(t),
$$  
that have zero mean and are delta correlated
\begin{eqnarray}
\label{eq:se}
\langle f_{Ej}(t)\rangle=\langle f_{Nj}(t)\rangle&=&0,\nonumber \\
\langle f_{Ej}^R(t)f_{Ej}^I(t)\rangle&=&0,\nonumber \\
\langle f_{Ei}^R(t)f_{Ej}^R(t')\rangle = \langle f_{Ei}^I(t)f_{Ej}^I(t')\rangle&=&D_E\delta_{ij}\delta(t-t'), \\
\langle f_{Ni}(t)f_{Nj}(t')\rangle&=&2D_N\delta_{ij}\delta(t-t').\nonumber
\end{eqnarray}
Here, $\delta_{ij}$ is the Kronecker delta and $\delta(t-t')$ is the Dirac delta function. 
In the calculations we use $D_E=0.05$ and $D_N=3.5\times10^{-8}$~\cite{WIE08}.

To measure the quality of synchronisation we introduce the  {\it order parameter}, $I_M(t)$, 
and the {\it average  order parameter}, $\langle I_M \rangle$, as
\begin{eqnarray}
\label{eq:intensity}
\langle I_M\rangle = \lim_{T\rightarrow\infty}\frac{1}{T}\int_0^T I_M(t)\,dt=
\lim_{T\rightarrow\infty}\frac{1}{T}\int_0^T\left| \sum_{j=1}^M E_j(t)\right|^2 dt. 
\end{eqnarray}
The physical meaning of $I_M(t)$ and $\langle I_M\rangle$ is illustrated in Fig.~\ref{fig:1b}.
If $M$ identical lasers are placed at an equal distance from a small (the order of a 
wavelength) spot and their light is focused onto this spot, 
then $I_M(t)$ and $\langle I_M\rangle$ are the instantaneous and average light intensity at the spot, respectively.
A single laser  oscillates with a random phase owing to spontaneous emission noise so  
that, for independent lasers, $\langle I_M\rangle$
is proportional to $M$ times the average intensity of a single laser. This follows 
directly from Eq.~(\ref{eq:intensity})
assuming lasers with identical amplitudes, $|E_j(t)|$, and uncorrelated random phases, 
$\arg(E_j(t))$. 
However, when the lasers oscillate in phase,  
one expects  $\langle I_M\rangle$ to be equal  $M^2$ times the average 
intensity of a single laser. This follows directly from Eq.~(\ref{eq:intensity})
assuming lasers with identical amplitudes and  phases.
We speak of synchronisation when $\langle I_M\rangle\approx M^2$, 
different degrees of partial synchronisation when $M<\langle I_M\rangle<M^2$, and lack of synchronisation
when $\langle I_M\rangle\approx M$. Note that $\langle I_M\rangle> M^2$ indicates 
trivial synchronisation, where the  external forcing term, $f_{ext}(t)$ becomes 'larger' than 
the oscillator terms on the right-hand side of Eq.~(\ref{eq:E2}). 
For comparability reasons, we now briefly review 
the case of a monochromatic forcing and then move on to the case of 
stochastic forcing.

\begin{figure}[t!]
\includegraphics[width=11cm]{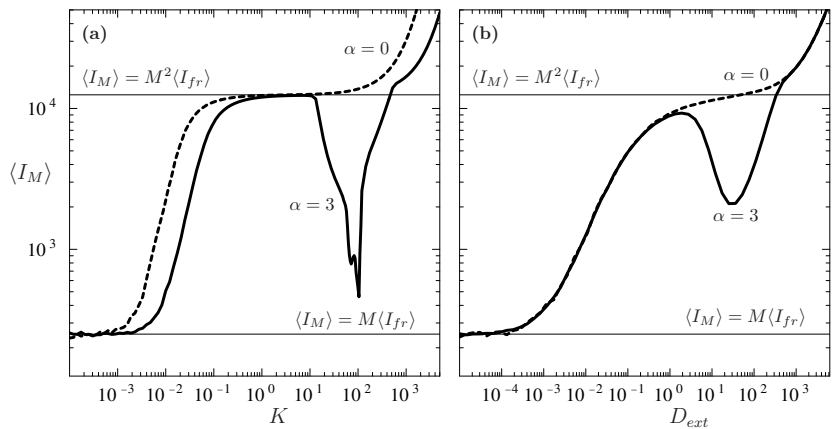}
\caption{
The average order parameter as defined by Eq.~(\protect\ref{eq:intensity}) for $M=50$ uncoupled
lasers with common (a) monochromatic and (b) white noise external forcing vs. 
the forcing strength for (dashed) $\alpha=0$
and (solid) $\alpha=3$; $J=5$ and $\nu_{ext}=\Delta$ in  panel (a).
Adapted from~\protect\cite{WIE09} with permission.
}
\label{fig:2}

\end{figure}

Let us consider a {\it monochromatic external forcing}
$$
f_{ext}(t)=Ke^{i\nu_{ext} t},
$$ 
where $K\in \mathbb{R}$ is the forcing strength and 
$\nu_{ext}$ is the detuning (difference) between 
the reference frequency chosen for $\Delta$ in Eq.~(\ref{eq:E}) and
the forcing frequency. 
Such an external forcing breaks the $\mathbb{S}$$^1$-symmetry and can force each 
laser to fluctuate in the vicinity of the well-defined external forcing 
phase, $\nu_{ext} t$, as opposed to a random walk. This phenomenon was studied 
in~\cite{CHO75} as a thermodynamic phase transition. 
Figure~\ref{fig:2}(a) shows $\langle I_{M}\rangle$ versus $K$, for an external 
forcing resonant with the laser,
$$
\nu_{ext}=\Delta.
$$ 
\begin{figure}[t!]
\includegraphics[width=9cm]{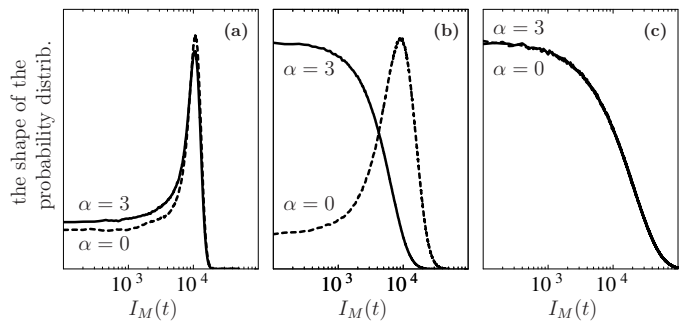}
\caption{
The shape of the probability distribution of $I_M(t)$ for  white noise external forcing;
(dashed) $\alpha=0$ and  (solid) $\alpha=3$. From (a) to (c) 
$D_{ext}= 10^0, 10^1$, and $10^3$. Adapted from~\protect\cite{WIE09} with permission.
}
\label{fig:3}
\end{figure}
\noindent
Because of the intrinsic spontaneous emission noise, the forcing amplitude has 
to reach a certain threshold before synchronisation occurs.
For $\alpha=0$, a sharp onset of synchronisation at $K\approx10^{-3}$ is 
followed by a wide range of $K$ with  synchronous behaviour where 
$\langle I_{M}\rangle= M^2 \langle I_{fr}\rangle$. Here,
$\langle I_{fr}\rangle$ is the average intensity of a single laser without forcing. 
At around $K=10^{2}$,  $\langle I_{M}\rangle$ starts increasing above 
$M^2 \langle I_{fr}\rangle$. Whereas lasers still remain synchronised, 
this increase indicates that the external forcing is no longer 'weak'. Rather, it 
becomes strong enough to cause an increase in the average intensity 
of each individual laser. A very different scenario is observed for $\alpha=3$. 
There, the onset of synchronisation is followed by an 
almost complete loss of synchrony just before $\langle I_{M}\rangle$ increases above 
$M^2 \langle I_{fr}\rangle$. The loss of synchrony is caused by  externally induced 
bifurcations and ensuing chaotic dynamics. These bifurcations have been studied in detail, 
both theoretically~\cite{WIE05,ERN96,CHL03,BON07} and experimentally~\cite{WIE02}, 
and are well understood.

The focus of this work is synchronisation to {\it white noise external forcing} represented 
by the  complex random process that is Gaussian, has zero mean, and is delta correlated
\begin{eqnarray}
\label{eq:wn}
f_{ext}(t)&=&f_{ext}^R(t)+if_{ext}^I(t),\nonumber \\
\langle f_{ext}(t)\rangle&=&\langle f_{ext}^R(t)f_{ext}^I(t)\rangle=0,\\
\langle f_{ext}^R(t)f_{ext}^R(t')\rangle &=& \langle f_{ext}^I(t)f_{ext}^I(t')\rangle =D_{ext}\delta(t-t')\nonumber.
\end{eqnarray}
White noise synchronisation is demonstrated in Fig.~\ref{fig:2}(b) 
where we plot $\langle I_{M}\rangle$ versus  $D_{ext}$. 
For $\alpha=0$,  a clear onset of synchronisation  at around $D_{ext}=10^{-3}$ is
followed by synchronous behaviour at larger $D_{ext}$. 
In particular, there exists a range of $D_{ext}$ where white noise external forcing 
is strong enough to synchronise phases of intrinsically  noisy lasers 
but weak enough so that each individual laser has small
intensity fluctuations and its average intensity remains unchanged.
In the probability distributions for  $I_M(t)$ in Fig.~\ref{fig:3}(a--b), the distinct peak
at $I_M(t)\approx M^2 \langle I_{fr}\rangle$
and a noticeable tail at smaller $I_M(t)$ indicate 
synchronisation that is not perfect. 
Rather, synchronous behaviour is occasionally interrupted with
short intervals of asynchronous behaviour owing to different intrinsic noise within each laser.
For $D_{ext}> 10^2$ the external forcing is no longer 'weak' and causes an increase in the
intensity fluctuations and the average intensity of each individual laser.
Although the lasers remain in synchrony, $\langle I_M\rangle$ increases 
above $M^2 \langle I_{fr}\rangle$ [Fig.~\ref{fig:2}(b)]
and exhibits large fluctuations [Fig.~\ref{fig:3}(c)] as in the asynchronous case.
A very different scenario is  observed again for $\alpha=3$. 
There, the onset of synchronisation is followed by a
significant loss of synchrony for $D_{ext}\in(4, 100)$. 
In this range of the forcing strength, one finds qualitatively different dynamics for $\alpha=0$ and $\alpha=3$ 
as revealed by different probability distributions in Fig.~\ref{fig:3} (b).

Interestingly, comparison between (a) and (b) in Fig.~\ref{fig:2}  shows that
some general aspects of synchronisation to a monochromatic and white noise external forcing
are strikingly similar. In both cases there is a clear onset of synchronisation followed by a significant 
loss of synchrony for sufficiently large $\alpha$, and subsequent revival of synchronisation for 
stronger external forcing. However, the dynamical mechanism responsible for the loss of synchrony 
in the case of white noise external forcing has not been fully explored.

\section{Definition of noise synchronisation}
\label{sec:sb}

The previous section motivates further research to reveal the dynamical  mechanism
responsible for the loss of synchrony observed in Fig.~\ref{fig:2}(b). To facilitate the analysis,
we  define synchronisation to irregular external forcing within the framework of 
random dynamical systems. Let us 
consider a $n$-dimensional, nonlinear, dissipative, autonomous dynamical system, referred to 
as {\it unforced system}
\begin{eqnarray}
\label{eq:osc}
\frac{dx}{dt}&= f(x,p),
\end{eqnarray}
where $x\in{\mathbb R}^n$ is the state vector and  $p\in{\mathbb R}^k$ is the parameter vector that 
does not change in time. An {\it external forcing} is denoted with  $f_{ext}(t)$, and 
the corresponding non-autonomous {\it forced system} reads
\begin{eqnarray}
\label{eq:fosc}
\frac{dx}{dt}&= f(x,p)+ f_{ext}(t).
\end{eqnarray}
Let $x(t,t_0,x_0)$ denote a {\it trajectory} or {\it solution} of~(\ref{eq:fosc}) that passes 
through $x_0$ at some initial time $t_0$. In situations where explicitly displaying the initial condition 
is not important we denote the trajectory simply as $x(t)$. For an infinitesimal displacement 
$\delta{x(0)}$ from $x(0,t_0,x_0)$, the {\it largest Lyapunov exponent} along $x(t,t_0,x_0)$ is given by
\begin{eqnarray}
\label{eq:le}
\lambda_{max}=\lim_{t\rightarrow\infty}\;\frac{1}{t}\,\ln\,\frac{|\delta{x(t)|}}{|\delta{x(0)}|}.
\end{eqnarray}
If the external forcing, $f_{ext}(t)$, is stochastic, Eq.~(\ref{eq:fosc}) defines 
a {\it random dynamical system} where $\lambda_{max}$ does not depend on  
the noise realisation, $f_{ext}$~\cite{ARN98}. Furthermore, we define:\\

\noindent
{\bf Definition 1.} An {\it (self-sustained) oscillator}  is an unforced system~(\ref{eq:osc}) 
with a  stable hyperbolic limit cycle.  
\\

\noindent
{\bf Definition 2.} An attractor for the forced system~(\ref{eq:fosc}) with 
stochastic forcing $f_{ext}(t)$ is called a {\it random sink} ({\bf rs}) if 
$\lambda_{max}<0$, and a {\it random strange attractor} ({\bf rsa}) 
if $\lambda_{max}>0$. 
\\

\noindent
{\bf Definition 3.} A {\it stochastic d-bifurcation} is a qualitative change in the random 
attractor when $\lambda_{max}$  crosses through zero~\cite[Ch.9]{ARN98}.\\

\noindent
{\bf Definition 4.} An oscillator is {\it synchronised} to stochastic 
forcing $f_{ext}(t)$ on a bounded subset $D\subset {\mathbb R}^n$ if the corresponding forced 
system~(\ref{eq:fosc}) has a random sink in the form of a unique attracting 
trajectory, $a(t,f_{ext})$, such that 
$$
\displaystyle \lim_{t_0\rightarrow-\infty}| x(t,t_0,x_0) - a(t,f_{ext})|\rightarrow 0,\\
$$ 
for fixed $t>t_0$ and all $x_0\in D$.
\\

\noindent
By Definition 1, an unforced oscillator has zero $\lambda_{max}$ on an open set of parameters.
In the presence of stochastic external forcing,  $\lambda_{max}$ becomes 
either positive or negative for typical parameter values~\cite{ARN98,GOL06} and 
 remains zero only at some special parameter values defining {\it stochastic d-bifurcations}.  
Synchronisation in  Definition 4 is closely related to 
{\it generalised synchronisation}~\cite{RUL95,ABA96} or {\it weak synchronisation}~\cite{PYR96} 
--- a phenomenon that requires a time-independent 
functional relationship between the measured properties of the forcing and the 
oscillator~\cite{BRO00}.
Following Refs.~\cite[Ch.9]{ARN98} and~\cite{LAN02},  
we used in Definition 4 the notion of {\it pullback convergence} where the asymptotic 
behaviour is studied for $t_0\rightarrow -\infty$ and fixed $t$. (For a study of different 
notions of convergence in  random dynamical systems we refer the reader to Ref.~\cite{ASH03}.) 
Whereas $\lambda_{max}$ does not depend on the noise realisation, $f_{ext}$, random sinks 
and random strange attractors do depend on $f_{ext}$. 
Hence the $f_{ext}$ dependence in $a(t,f_{ext})$ in Definition 4.
Since $\lambda_{max}<0$ does not imply a unique attracting trajectory it is 
not sufficient to show synchronisation as defined in Definition 4. 
In general, there can be a number of coexisting (locally) attracting trajectories that 
belong to a {\it global pullback attractor}~\cite{LAN02}. In such cases, one can choose $D$ to lie in 
the basin of attraction of one of the locally attracting trajectory and speak of synchronisation on $D$.




\begin{figure}[h!]
\includegraphics[width=8.4cm]{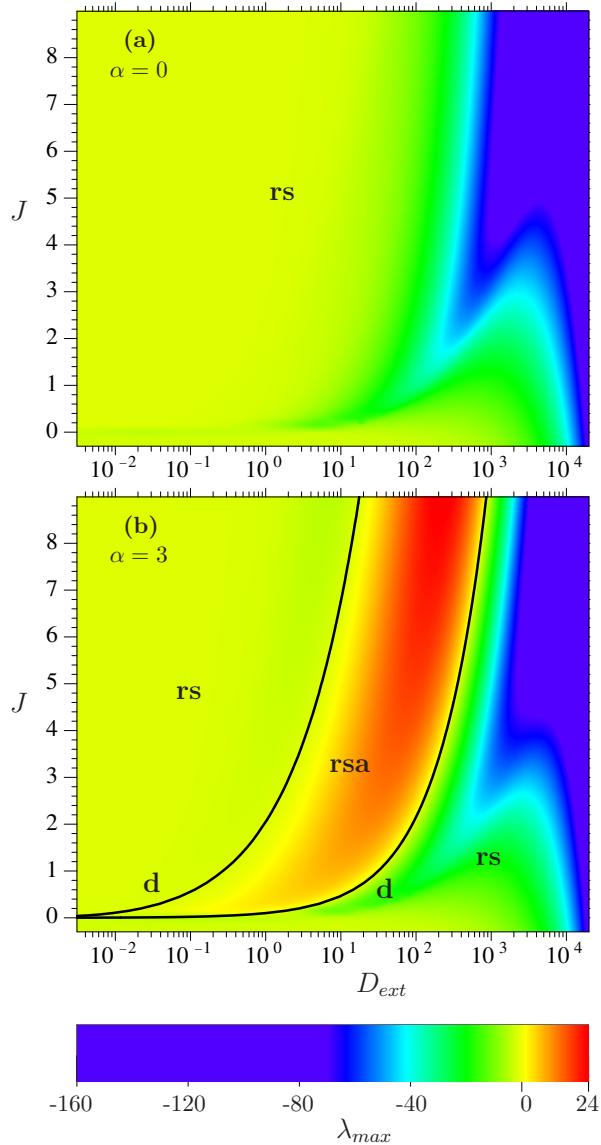}
\caption{
Bifurcation diagram for the white noise forced laser model~(\protect\ref{eq:E3}--\protect\ref{eq:N3}) in the 
$(D_{ext},J)$-plane for (a) $\alpha=0$ and (b) $\alpha=3$. 
{\bf rs} stands for random sink, {\bf rsa} stands for random strange attractor, the two black curves 
{\bf d} denote stochastic d-bifurcation, and colour coding is for $\lambda_{max}$.
An unforced laser has a Hopf bifurcation at  $(D_{ext}=0,J=0)$. Adapted from~\protect\cite{WIE09} with permission.
}
\label{fig:4}
\end{figure}

\begin{figure}[h!]
\includegraphics[width=11cm]{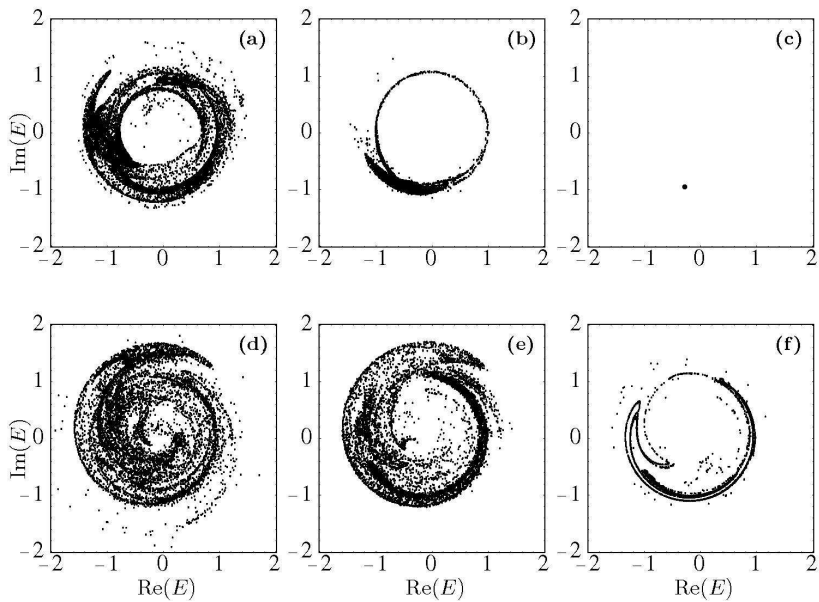}
\caption{
Pullback convergence to (a--c) a random sink for $D_{ext}=0.1$ and 
(d--f) a random strange attractor for $D_{ext}=0.5$ in the
white noise forced laser model~(\protect\ref{eq:E3}--\protect\ref{eq:N3}) 
with $\alpha=3$ and $J=1$. Shown are snapshots 
of $10\,000$ trajectories in projection onto the complex $E$-plane at time $t=30$. 
The initial conditions are uniformly distributed 
on the $E$-plane with $N=0$ at different initial times $t_0$=(a,d) 29, (b,e) 28, and (c,f) 0.
}
\label{fig:5}
\end{figure}

\section{Synchronisation transitions via Stochastic d-bifurcation}
\label{sec:bifdiag}

To facilitate the analysis we make use of Definitions 1--3 and, henceforth, consider noise 
synchronisation in the laser model with white noise external forcing
\begin{eqnarray}
\label{eq:E3}
\frac{dE}{dt}&=& i\Delta\,E+ g\gamma(1-i\alpha)NE +f_{ext}(t),\\
\label{eq:N3}
\frac{dN}{dt}&=& J - N - (1+g N)|E|^2 ,
\end{eqnarray}
but without the intrinsic spontaneous emission noise.  
Now, owing to the absence of intrinsic noise, Definition 4 is equivalent to the 
synchronisation detection scheme chosen in 
Section~\ref{sec:synch}. More specifically, the evolution of $M$ trajectories 
starting at different initial conditions for a single laser with external forcing is the same as 
the evolution of an ensemble of $M$ identical  uncoupled lasers with the same forcing, 
where each laser starts at a different initial condition. 
A random sink for~(\ref{eq:E3}--\ref{eq:N3})
in the form of a unique attracting trajectory, $a(t,f_{ext})$, makes trajectories for 
different initial conditions converge to each other. In an ensemble of $M$ identical 
uncoupled lasers with common forcing this means that
$\langle I_{M}\rangle\rightarrow  M^2 \langle I_{fr}\rangle$ in time
so that synchronisation is detected.  
A random strange attractor for~(\ref{eq:E3}--\ref{eq:N3}), where nearby trajectories separate 
exponentially fast because $\lambda_{max}>0$, implies $\langle I_{M}\rangle<  M^2 \langle I_{fr}\rangle$ so that 
incomplete synchronisation or lack of synchronisation is detected.

Figure~\ref{fig:4} shows effects of white noise external forcing on the 
sign of the otherwise zero 
$\lambda_{max}$ in an unforced laser. For $\alpha=0$, external forcing always 
shifts  $\lambda_{max}$
to negative values meaning that the system has a random sink for $D_{ext}>0$ 
and $J>0$ [Fig.~\ref{fig:4}(a)]. Additionally, this random sink is a unique trajectory, $a(t,f_{ext})$, 
meaning that the laser is synchronised to white noise external forcing.
However, for $\alpha=3$ there are two curves of stochastic 
d-bifurcation where $\lambda_{max}$ crosses through zero [Fig.~\ref{fig:4}(b)]. 
Noise synchronisation is lost for parameter settings between these two curves, where $\lambda_{max}>0$
indicates a random strange attractor .
Pullback convergence to two qualitatively different random attractors 
found  for $\alpha=3$ is shown in Fig.~\ref{fig:5}. At fixed time $t=30$, we
take snapshots of $10\,000$ trajectories for a grid of initial conditions
with different initial times $t_0$.
In (a--c), trajectories converge in the pullback sense to a random sink. 
The random sink appears in the snapshots as a single dot whose position is different for different $t$ 
or different noise realisations $f_{ext}$. 
In (d--f), trajectories converge  in the pullback sense to a random strange 
attractor that appears in the snapshots as a fractal-like structure.  This structure 
remains fractal-like but is different for different $t$ or different noise realisations $f_{ext}$.

\subsection{Class-B Laser Model vs. Landau-Stuart Equations}
\label{sec:bifdiagcompar}
The stochastic d-bifurcation uncovered in the previous section  has been reported  in 
biological systems~\cite{KOS04,GOL05,GOL06,LIN07,LIN08} and should appear in a general class
of  oscillators with stochastic forcing. 
Here, we use the laser model in conjunction with the Landau-Stuart model
to  address its dependence on the three  parameters: $D_{ext},J,$ 
and $\alpha$, and to uncover its universal properties. With an exception of certain 
approximations~\cite{GOL06}, 
this problem is beyond the reach of analytical techniques and so 
numerical analysis is the tool of choice.

\begin{figure}[t!]
\includegraphics[width=11.5cm]{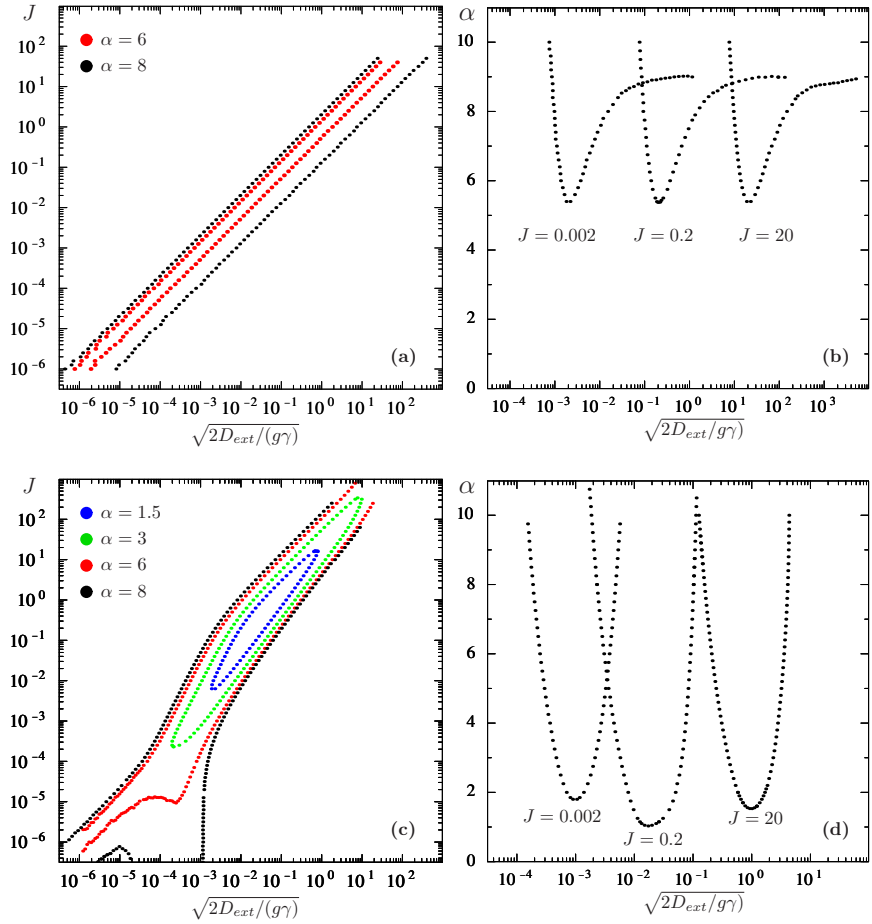}
\caption{
A three-parameter study of stochastic d-bifurcation (a--b) in the  
white noise forced Landau-Stuart model~(\protect\ref{eq:reducedn}) and (c--d) the 
white noise forced laser model~(\protect\ref{eq:E3}--\protect\ref{eq:N3}). To facilitate the 
comparison between~(\protect\ref{eq:reducedn}) and~(\protect\ref{eq:E3}--\protect\ref{eq:N3}), we plotted  results 
for~(\ref{eq:E3}--\ref{eq:N3}) with the rescaled forcing strength, $\sqrt{2D_{ext}/g\gamma}$. 
Adapted from~\protect\cite{WIE09} with permission.
}
\label{fig:6}
\end{figure}

To help identifying  effects characteristic to the more complicated laser model, we first consider the 
Landau-Stuart model with white noise external forcing
\begin{eqnarray}
\label{eq:reducedn}
\frac{dE}{d\tilde{t}}= 
\left[J +i\left(\tilde{\Delta}-\alpha(J -|E|^2) \right) \right]E - E|E|^2 +f_{ext}(\tilde{t}),
\end{eqnarray}
where, for the rescaled time $\tilde{t}$, the external forcing correlations become
$$
\langle f_{ext}^R(\tilde{t})f_{ext}^R(\tilde{t}')\rangle = 
\langle f_{ext}^I(\tilde{t})f_{ext}^I(\tilde{t}')\rangle =
\frac{D_{ext}}{g\gamma}\delta(\tilde{t}-\tilde{t}')\nonumber.
$$
Figure~\ref{fig:6}(a--b) shows the dependence of the 
d-bifurcation on $J$, $D_{ext}$, and $\alpha$ in Eq.~(\ref{eq:reducedn}). In the three-dimensional 
$(J, D_{ext},\alpha)$-parameter space, the two-dimensional surface of d-bifurcation appears to originate from the half 
line $(D_{ext}=0,J=0,\alpha> 5.3)$ of 
the deterministic Hopf bifurcation, has a ridge at 
$\alpha_{min}\approx 5.3$,  and is asymptotic to $\alpha\approx 9$
with increasing $D_{ext}$ [Fig.~\ref{fig:6}(b)]. 
Furthermore, numerical results in Fig.~\ref{fig:6}(b) suggest that the shape of the 
d-bifurcation curve in the two-dimensional section $(D_{ext},\alpha)$ is independent of $J$.
As a consequence, for fixed $\alpha$ within the range $\alpha\in(5.3,9)$ 
one finds two d-bifurcation curves in the $(D_{ext},J)$-plane [Fig.~\ref{fig:6}(a)] that are parametrised by   
\begin{eqnarray}
\label{eq:scalling}
J_j= C_{j}(\alpha)\sqrt{2D_{ext}}\,\,,\;\;\;\;\mbox{where}\;\;\;\;j=1,2, 
\end{eqnarray}
and bound the region with a random strange attractor.
Since $C_{1}(\alpha_{min})=C_{2}(\alpha_{min})=1$, these two curves merge into 
a single curve
\begin{eqnarray}
\label{eq:scalling2}
J= \sqrt{2D_{ext}}\,\,,
\end{eqnarray}
when $\alpha=\alpha_{min}$.
On the one hand, for  $\alpha\le\alpha_{min}$, the region with a random strange attractor 
disappears from the $(D_{ext},J)$-plane.
On the other hand, for $ \alpha> 9$, there is just one  d-bifurcation curve in 
the $(D_{ext},J)$-plane, meaning that the region with a random strange attractor 
becomes unbounded towards increasing $D_{ext}$ [Fig.~\ref{fig:6}(b)].

Similar results are expected for any white noise forced  oscillator with shear that is  near 
a Hopf bifurcation, and for 'weak' forcing. This claim is supported with numerical analysis of the laser 
model~(\ref{eq:E3}--\ref{eq:N3}) in Fig.~\ref{fig:6}(c--d). 
For a fixed $\alpha$, 
Eq.~(\ref{eq:reducedn}) and Eqs.~(\ref{eq:E3}--\ref{eq:N3}) 
give identical results if 
the forcing  is weak enough but significant discrepancies 
arise with increasing forcing strength.
First of all, it is possible to have one-dimensional sections of the 
$(D_{ext},J)$-plane for fixed $J$
 with multiple uplifts of $\lambda_{max}$ to positive values 
[black dots for $J<10^{-6}$ in Fig.~\ref{fig:6}(c)].
Secondly, the parameter region with a random strange attractor for  Eqs.~(\ref{eq:E3}--\ref{eq:N3})  
expands towards much lower values of $\alpha>\alpha_{min}\approx1$  [Fig.~\ref{fig:6}(d)].
 Thirdly,  the shape of the two-dimensional 
surface of d-bifurcation in the laser model becomes dependent on $J$ and
has a minimum rather than a ridge. 
As a consequence, 
although the stochastic d-bifurcation 
seems to originate from the half line $(D_{ext}=0,J=0,\alpha> 5.3)$, 
it will appear in the $(D_{ext},J)$-plane even
for $\alpha\in(1,5.3)$ as a closed and isolated curve away from the origin of 
this plane [Fig.~\ref{fig:6}(c)]. Finally, the  region of  random strange attractor
remains bounded  in the $(D_{ext},J)$-plane even for large $\alpha$.

\begin{figure}[t!]
\includegraphics[width=7.3cm]{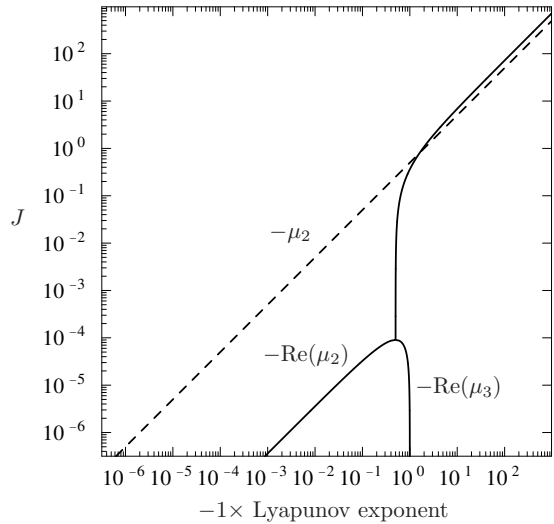}
\caption{
$J$ versus (solid) the two non-zero Lyapunov exponents, $\mbox{Re}(\mu_2)$ and $\mbox{Re}(\mu_3)$, 
for the limit cycle in the laser model~(\protect\ref{eq:E}--\protect\ref{eq:N}),
and (dashed) the nonzero Lyapunov exponent, $\mu_2$, for the limit cycle in the Landau-Stuart model~(\protect\ref{eq:reduced}). Adapted from~\protect\cite{WIE09} with permission.
}
\label{fig:7}
\end{figure}

To unveil the link between the transient dynamics of unforced systems and  the forcing-induced stochastic d-bifurcation,
we plot  $J$ versus Lyapunov 
exponents in Fig.~\ref{fig:7}; 
note that  Lyapunov exponents, $\lambda_i$, and Floquet exponents, $\mu_i$, are 
related by $\lambda_i={\textrm Re}[\mu_i]$.
A comparison between  Figs.~\ref{fig:6} and \ref{fig:7} 
shows strong correlation between  the relaxation towards the limit cycle 
and the d-bifurcation. 
In the Landau-Stuart model~(\ref{eq:reduced}), 
the linear relation~(\ref{eq:fels}) between  $J$  and  
the non-zero Lyapunov exponent, $\mu_2$ (dashed line in Fig.~\ref{fig:7}),
results in a linear parametrisation~(\ref{eq:scalling}) of d-bifurcation curves 
in the ($J, \sqrt{2D_{ext}})$-plane [Fig.~\ref{fig:6}(a)].
In the laser model~(\ref{eq:E}--\ref{eq:N}), the nonlinear 
relation~(\ref{eq:mu1}--\ref{eq:mu2}) between $J$ and the  non-zero 
Lyapunov exponents, $\mbox{Re}(\mu_{2})$ and $\mbox{Re}(\mu_{3})$ 
(solid curves in Fig.~\ref{fig:7}), results in a very 
similar nonlinear parametrisation of d-bifurcation curves in the 
($J, \sqrt{2D_{ext}}$)-plane [Fig.~\ref{fig:6}(c)]. The
splitting up of the chaotic region bounded by the black dots for $J<9\times10^5$ 
in Fig.~\ref{fig:6}(c)
is related to two different  eigendirections normal to the limit cycle
with significantly different timescales of transient dynamics towards the limit cycle
(the two corresponding Lyapunov exponents are shown in solid in Fig.~\ref{fig:7}). 
Finally, the appearance of relaxation oscillations in the laser system is associated with a 
noticeable expansion of the chaotic region, in particular,  towards smaller $\alpha$.



\section{ Noise-induced strange attractors}
\label{sec:interpret}

Complicated invariant sets,  such as strange attractors, 
require a balanced interplay between phase space expansion and contraction~\cite{GUC83}. 
If phase space expansion in certain directions is properly compensated 
by phase space contraction in some other directions, nearby trajectories can
separate exponentially fast ($\lambda_{max}>0$) and yet remain within a bounded 
subset of the  phase space.

\begin{figure}[t!]
\includegraphics[width=9cm]{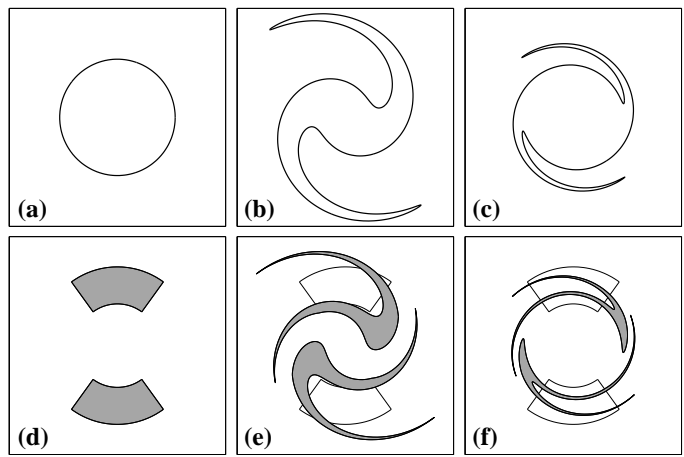}
\caption{
 Time evolution of sets of initial conditions showing the
creation of horseshoes in the phase space of a suitably kicked
laser model with no shear ($\alpha=0$). The sets of initial conditions are (a) the
stable circle and (d) boxes containing parts of the circle.
Shown are phase portraits  [(a) and (d)] before,  [(b) and (e)] immediately after the first kick, and [(c) and (f)] some time after the first kick.
}
\label{fig:8}
\end{figure}

It has been recently proven that, when
{\it suitably} perturbed, any stable hyperbolic limit cycle can be
turned into `observable' chaos (a strange attractor)~\cite{WAN03}. This result is
derived for periodic discrete-time perturbations (kicks) that
deform the stable limit cycle of the unkicked system. The key concept
is the creation of horseshoes via a stretch-and-fold action due to an
interplay between the kicks and the local geometry of the phase
space. Depending on the degree of shear, quite different kicks are
required to create a stretch-and-fold action and
horseshoes. Intuitively, it can be described as follows. In systems
without shear, where points in phase space rotate with the same
angular frequency about the origin of the complex $E$-plane 
independent of their distance from the origin, the
kick alone has to create the stretch-and-fold action. This is
demonstrated in Fig.~\ref{fig:8}. Horseshoes are formed 
as the system is suitably kicked in both radial and angular directions
 and then relaxes back to the attractor [the circle in Fig.~\ref{fig:8}(a)] of the
unkicked system. Repeating this process reveals chaotic
invariant sets. However, showing rigorously that a specific kick results in
`observable' chaos is a non-trivial task~\cite{WAN03}. 
In the presence of shear, where points
in phase space rotate with different angular frequencies depending on
their distance from the origin, the kick does not have to be so
specific or carefully chosen. In fact, it may be sufficient to kick
non-uniformly in the radial direction alone, and rely on natural
forces of shear to provide the necessary stretch-and-fold action.

\begin{figure}[t!]
\includegraphics[width=11cm]{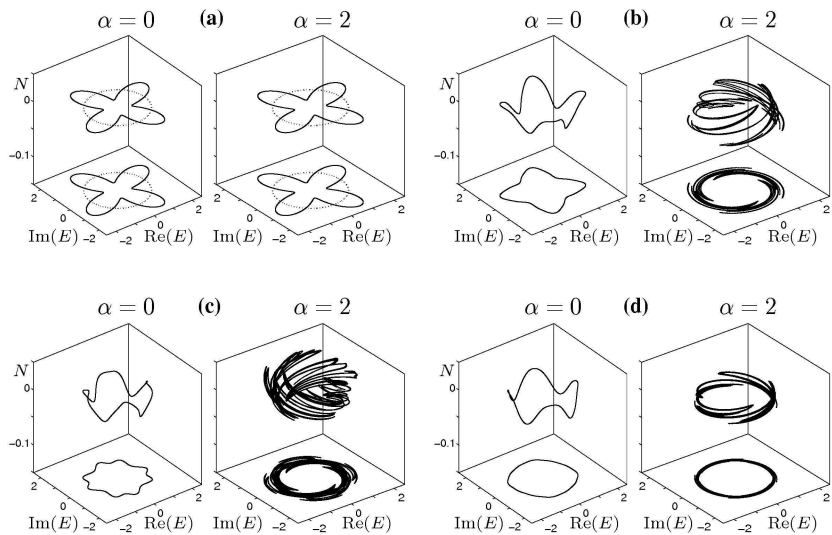}
\caption{
Snapshots at times (a) $t=0$, (b) $t=0.35$, (c) $t=0.8$,
and (d) $t=1$ showing the time evolution of $15 000$ trajectories with
initial conditions distributed equally over the stable 
circle for~(\protect\ref{eq:E}--\protect\ref{eq:N}). Kicks in the radial direction alone are applied at
times $t=0,0.25,0.5,0.75$. A comparison between $\alpha =0$ and $\alpha =2$
illustrates the $\alpha$-induced stretch-and-fold action in the laser
phase space.
}
\label{fig:9}
\end{figure}

These effects are illustrated in
Fig.~\ref{fig:9} for the single laser
model~(\ref{eq:E}--\ref{eq:N}) with non-uniform kicks in the radial
direction alone for $\alpha=0$ (no shear) and $\alpha=2$ (shear). There, we set $\Delta=0$
and refer to the stable limit cycle [dashed circle in Fig.~\ref{fig:9}(a)] as  $\Gamma$. 
Kicks modify the electric field amplitude, 
$|E|$, by a
factor of $0.8\sin[4\,\arg(E)]$ at times $t = 0,\, 0.25,\,
0.5,\,\mbox{and}\, 0.75$, but leave the phase, $\arg(E)$, unchanged. 
For $\alpha=0$ each point on the black curve
spirals onto $\Gamma$ in time but remains within the same isochrone
defined by a constant electric-field phase, $\arg(E)=\arg(E(0))$. Hence, the black
curve does not have any folds at any time. However, for $\alpha =2$,
a kick moves most points on the black curve to  different isochrones so that
 points with larger amplitudes $|E|$ rotate with
larger angular frequencies. This gives rise to an intricate
stretch-and-fold action. Folds and
horseshoes can be formed under the evolution of the flow even though the
kicks are in the radial direction alone.

\begin{figure}[t!]
\includegraphics[width=7.5cm]{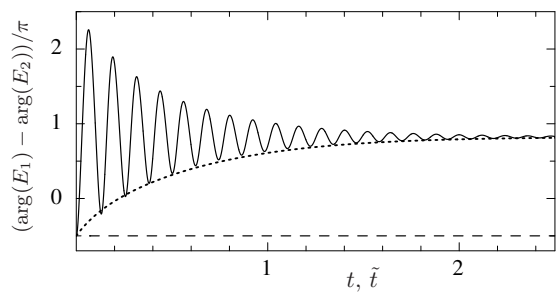}
\caption{
 Phase-space stretching along the limit cycle shown as time evolution of the phase difference  
between two trajectories starting at different isochrones.
The three curves are obtained using (dashed) Eqs.~(\protect\ref{eq:E}--\protect\ref{eq:N}) with 
$\alpha=0$, (solid) Eqs.~(\protect\ref{eq:E}--\protect\ref{eq:N}) with $\alpha=3$, and 
(dotted) Eq.~(\protect\ref{eq:reduced}) with $\alpha=3$. Adapted from~\protect\cite{WIE09} with permission. 
}
\label{fig:10}
\end{figure}

In the laser model, stretch-and-fold action is 
significantly enhanced by the spiralling transient motion about $\Gamma$. 
For $J> 10^{-1}$, the laser model~(\ref{eq:E}--\ref{eq:N}) and the Landau-Stuart 
model~(\ref{eq:reduced}) have nearly identical relaxation timescales toward $\Gamma$ 
(Fig.~\ref{fig:7}).  
However,  owing to one additional degree of freedom and 
oscillatory relaxation~(\ref{eq:mu2}),  the instantaneous stretching along $\Gamma$ in the 
three-dimensional laser vector field~(\ref{eq:E}--\ref{eq:N}) can be much stronger compared 
to the planar vector field~(\ref{eq:reduced}), especially at short times after the perturbation.
This effect is illustrated in Fig.~\ref{fig:10} by the time evolution of the phase difference, 
$\arg(E_1(t))-\arg(E_2(t))$, between two trajectories, 1 and 2, starting at different 
isochrones for $\alpha=3$. 
Since both vector fields have identically shaped isochrones~(\ref{eq:isochrone3}), 
the phase difference converges to the 
same value as time tends to infinity.  However, at small $t$, the oscillatory phase 
difference for~(\ref{eq:E}--\ref{eq:N}) can exceed significantly the monotonically 
increasing phase difference for~(\ref{eq:reduced}) (compare  solid and dotted 
curves in Fig.~\ref{fig:10}).

It is important to note that the rigorous results for turning stable
limit cycles into chaotic attractors are derived for periodic
discrete-time perturbations. Stochastic forcing is a
continuous-time perturbation, meaning that
the analysis in~\cite{WAN03} cannot be directly applied to our
problem. Nonetheless, such analysis gives a  valuable insight as to why random 
chaotic attractors appear for $\alpha$ sufficiently large, and it helps distinguish 
effects of stochastic forcing.

Here, we demonstrated that purely additive white noise forcing
is sufficient to induce random strange attractors in limit cycle oscillators.
Furthermore, numerical analysis in Sec.~\ref{sec:bifdiagcompar} shows that, 
in the case of stochastic forcing,
creation of strange attractors requires a different balance between 
the amount of shear, relaxation rate toward the limit cycle, 
and forcing strength, as compared to periodic forcing. Unlike in the case of 
discrete-time periodic forcing, the shear has to be strong enough, $|\alpha|>C>0$,
to allow sufficient stretch-and-fold action.
Provided that the shear is strong enough, 
the stochastic forcing strength needs to be at least comparable to the relaxation rate 
toward the limit cycle to allow formation of random strange attractors. 
Furthermore, we demonstrated that in higher dimensional systems, different eigendirections 
with distinctly different relaxation rates toward the limit cycle could give rise to more than
one region in the ($J, \sqrt{2D_{ext}}$)-plane with a random strange attractor. 
Last but not least, 
we revealed that the enhancement in the instantaneous stretch-and-fold action arising from 
laser relaxation oscillations results in a larger parameter region with a random strange attractor.

\section{Conclusions}
\label{sec:bif}

We used the class-B laser model in conjunction with
the Landau-Stuart model (Hopf normal form with shear) to study noise 
synchronisation and loss of synchrony via shear-induced stochastic d-bifurcations. 

We defined noise synchronisation in terms of pullback convergence of random attractors
and showed that a nonlinear oscillator can synchronise to stochastic external forcing. However, 
the parameter region with synchronous 
dynamics becomes interrupted with a single or multiple intervals of 
asynchronous dynamics if  amplitude-phase coupling or shear is 
sufficiently large. Stability analysis shows that the synchronous solution represented 
by a random sink 
loses stability via stochastic d-bifurcation to a random strange attractor. 
We performed a systematic study of this bifurcation with dependence on the three parameters: 
the Hopf bifurcation parameter (laser pump),  the amount of shear 
(laser linewidth enhancement factor), and the stochastic forcing strength.

In this way, we uncovered a vast 
parameter region with random strange attractors that are induced purely by stochastic forcing.
More specifically, in the three-dimensional parameter space, the two-dimensional surface of 
the stochastic bifurcation originates from the half-line of the deterministic Hopf bifurcation. 
In the plane of stochastic forcing strength and Hopf bifurcation parameter, one finds stochastic d-bifurcation curve(s)
bounding region(s) of random strange attractors if the amount of shear is sufficiently large. 
The shape of d-bifurcation curves 
is determined by the type and rate of the relaxation toward the limit cycle in the unforced oscillator.
Near the Hopf bifurcation 
and provided that stochastic forcing is weak enough,
the d-bifurcation curves satisfy the numerically uncovered power law~(\ref{eq:scalling}).
However, as the stochastic forcing strength increases, there might be deviations from this law. 
The deviations arise because different 
oscillators experience different effects of higher-order terms and additional degrees of freedom on the 
relaxation toward the limit cycle. 
In the laser example, the d-bifurcation curves 
deviate from the 
simple power law~(\ref{eq:scalling}) 
so that the region of a random strange attractor  splits up  and 
expands toward smaller values of shear as the forcing strength increases. 
We intuitively 
explained these results by demonstrating that the shear-induced stretch-and-fold action in 
the oscillator's phase space
facilitates creation of  horseshoes and strange attractors in response to external forcing. 
Furthermore, we showed that the  
stretch-and-fold action can be greatly enhanced by damped relaxation oscillations in the 
laser model, causing the 
deviation from the simple power law. 




\begin{thebibliography}{00}


\bibitem{PIK01} A.S. Pikovsky, M. Rosenblum, and J. Kurths,
{\it Synchronisation-a unified approach to nonlinear science},
Cambridge University Press, Cambridge, 2001.

\bibitem{PIK84} A.S. Pikovsky,in
{\it Nonlinear and turbulent processes in physics},
edited by R.Z. Sagdeev (Harwood Academic, Singapore, 1984),
1601--1604.

\bibitem{MAI95} Z.F. Mainen and T.J. Sejnowski,
``Reliability of spike timing in neocortical neurons'',
 {\em Science}  {\bf 268}  1503--1506 (1995).

\bibitem{JEN98} R.V. Jensen,
``Synchronization of randomly driven nonlinear oscillators'',
 {\em Phys. Rev. Lett.}  {\bf 58} R6907--R6910 (1998).

\bibitem{ZHO03} C.S. Zhou, J. Kurths, E. Allaria, S. Boccaletti, R. Meucci, and F. T. Arecchi,
``Constructive effects of noise in homoclinic chaotic systems'',
 {\em Phys. Rev. Lett.}  {\bf 67}  066220 (2003).

\bibitem{JUN04} J. Teramae and D. Tanaka,
``Robustness of the Noise-Induced Phase Synchronization in a General Class of Limit Cycle Oscillators'',
 {\em Phys. Rev. Lett.}  {\bf 93} 204103 (2004).

\bibitem{UCH04} A. Uchida, R. McAllister, and R. Roy,
``Consistency of nonlinear system response to complex drive signals'',
 {\em Phys. Rev. Lett.}  {\bf 93} 244102 (2004).


\bibitem{KOS04} E. Kosmidis and K. Pakdaman,
`` An Analysis of the Reliability Phenomenon in the FitzHugh-Nagumo Model'',
 {\em J. Comp. Neurosci.}  {\bf 14} 5--22 (2003).


\bibitem{GOL05} D.S. Goldobin and A. Pikovsky,
``Synchronisation and desynchronisation of self-sustained oscillators by common noise'',
 {\em Phys. Rev. E}  {\bf 71} 045201(R) (2005).

\bibitem{GOL06} D.S. Goldobin and A. Pikovsky,
``Antireliability of noise-driven neurons'',
 {\em Phys. Rev. E}  {\bf 73} 061906 (2006).


\bibitem{RUL95} N.F. Rulkov, M.M. Sushchik, L.S. Tsimring, and H.D.I. Abarbanel,
``Generalised synchronisation of chaos in directionally coupled chaotic systems'',
 {\em Phys. Rev. E}  {\bf 51} 980--994 (1995).


\bibitem{WIE09} S. Wieczorek
``Stochastic bifurcation in noise-driven lasers and Hopf oscillators'',
 {\em Phys. Rev. E}  {\bf 79} 036209 (2009).


\bibitem{LIN07}K.K. Lin, E. Shea-Brown, and L.-S. Young,
``Reliability of coupled oscillators'',
 {\em J. Nonlin. Sci. }  {\bf 19} 497--545 (2009).

\bibitem{LIN08} K.K. Lin, and L.-S. Young,
``Shear-induced chaos'',
 {\em Nonlinearity}  {\bf 21} 899--922 (2008).



\bibitem{WIE05} S. Wieczorek, B. Krauskopf, T.B. Simpson, and D. Lenstra,
``The dynamical complexity of optically injected semiconductor lasers'',
 {\em Phys. Rep.}  {\bf 416} 1--128 (2005).


\bibitem{WIE08} S. Wieczorek and W.W. Chow,
``Bifurcations and chaos in a semiconductor laser with coherent or noisy optical injection'',
 {\em Optics Communications} {\bf 282} (12) 2367--2379 (2009).


\bibitem{GUC83} J. Guckenheimer and P. Holmes,
{\em Nonlinear oscillations, dynamical systems, and bifurcations of vector fields},
Springer-Verlag, New York 1983.


\bibitem{HEN83} C.H. Henry,
``Theory of the linewidth of semiconductor laser'', 
{\em IEEE J. Quantum Electron.} {\bf QE--18} 259 (1982).


\bibitem{GUC74}  J. Guckenheimer,
``Isochrons and phaseless sets'',
 {\em J. Math. Biol.} {\bf 3} 259--273 (1974/75).


\bibitem{LAN82} R. Lang, 
``Injection locking properties of a semiconductor laser'',
{\it IEEE J. Quantum Electron.} {\bf 18}, 976 (1982).


\bibitem{CHO75} W.W. Chow, M.O. Scully, and E.W. van Stryland,
``Line narrowing in a symmetry broken laser'',
 {\em Opt. Commun.}  {\bf 15} 6--9 (1975).


\bibitem{ERN96} T. Erneux, V. Kovanis, A. Gavrielides, and P.M. Alsing,
``Mechanism for period-doubling bifurcation in a semiconductor laser subject to optical injection'',
 {\em Phys. Rev. A}  {\bf 53} 4372--4380 (1996).

\bibitem{CHL03} K.E. Chlouverakis and M.J. Adams,
``Stability map of injection-locked laser diodes using the largest Lyapunov exponent'',
 {\em Opt. Comm.}  {\bf 216} 405--412 (2003).

\bibitem{BON07} C. Bonatto and J.A.C. Gallas,
``Accumulation horizons and period adding in optically injected semiconductor lasers'',
 {\em Phys. Rev. E}  {\bf 75} 055204R (2007).

\bibitem{WIE02} S. Wieczorek, T.B. Simpson, B. Krauskopf, and D. Lenstra
``Global quantitative predictions of complex laser dynamics'',
 {\em Phys. Rev. E}  {\bf 65} 045207R (2002).

\bibitem{ARN98}
L. Arnold,
{\it Random Dynamical Systems},
(Springer-Verlag, Berlin Heidelberg, 1998).


\bibitem{ASH03}  P. Ashwin and G. Ochs,
``Convergence to local random attractors'',
 {\em Dyn. Syst.} {\bf 18} (2) 139--158 (2003).


\bibitem{RUL95} N.F. Rulkov, M.M. Sushchik,  L.S. Tsimring, and H.D.I. Abarbanel,
``Generalised synchronisation of chaos in directionally coupled chaotic systems'',
 {\em Phys. Rev. E} {\bf 51} (2) 980--994 (1995).

\bibitem{ABA96} H.D.I. Abarbanel N.F. Rulkov, and M.M. Sushchik,
``Generalised synchronisation of chaos: The auxiliary system approach'',
 {\em Phys. Rev. E} {\bf 53} (5) 4528--4535 (1996).

\bibitem{PYR96} K. Pyragas,
``Weak and strong synchronisation of chaos'',
 {\em Phys. Rev. E} {\bf 54} (5) R4508--R4511 (1996).

\bibitem{BRO00} R. Brown and L. Kocarev,
``A unifying definition of synchronisation for dynamical systems'',
 {\em Chaos} {\bf 10} (2) 344--349 (2000).


\bibitem{LAN02} J.A. Langa, J.C. Robinson, and A. S\'uarez,
``Stability, instability, and bifurcation phenomena in non-autonomous differential equations'',
 {\em Nonlinearity} {\bf 15} (2) 1--17 (2002).


\bibitem{WAN03} Q. Wang and L.-S. Young,
``Strange attractors in periodically-kicked limit cycles and Hopf bifurcations'',
 {\em Comm. Math. Phys.}  {\bf 240} 509--529 (2003).




\end{thebibliography}
\end{document}